\documentclass[preprint,12pt]{elsarticle}

\def\myrefs#1#2{ 
{\bigskip \noindent
{\Large \bf #2}  
 \list {[\arabic{enumi}]}{\settowidth\labelwidth{[#1]}
 \leftmargin\labelwidth 
 \advance\leftmargin\labelsep
 \usecounter{enumi} }  
 \def\newblock{\hskip .11em plus .33em minus .07em}
 \sloppy\clubpenalty4000\widowpenalty4000
 \sfcode`\.=1000\relax}  }

\def\nref#1{(\ref{#1})}

\def\newpage{\vfill\eject}
\def\comb#1,#2,{ \left( {#1 \atop #2 } \right)  }%
\def\prodd#1,#2,#3,{ \prod_{\scriptstyle #1 \atop\scriptstyle #2 }^{#3} }%
\def\summ#1,#2,#3,{ \sum_{\scriptstyle #1 \atop\scriptstyle #2 }^{#3} }%
\def\CC{\mathbb{C}}

\newtheorem{algor}{{\sc Algorithm}}[section]
\newtheorem{tabl}{Table}[section]

\def\betab{\begin{tabbing} 
xxxx\=xxxx\=xxx\=xx\=xx\=xx\=xx\=xx\=xx\=xx\=xx\=xx\=xx\= \kill} 
\def\entab{\end{tabbing}\vspace{-0.12in}}

\newcommand{\eq}[1]{\begin{equation}\label{#1}}
\newcommand{\en}{\end{equation}}

\newcommand{\beeq}[1]{\begin{equation}\label{#1}}
\newcommand{\eneq}{\end{equation}}

\journal{Computer Physics Communications}

\usepackage{amsmath,amssymb}
\usepackage{graphicx}
\usepackage{float}
\usepackage{algorithmic,algorithm}
\usepackage{bookmark}
\usepackage{url}
\usepackage{multirow}
\usepackage{hyperref}

\usepackage{amsthm}

\newcommand{\ab}{[\alpha,\beta]}

\newcommand{\bR}[1]{\mathbb{R}^{#1}}

\newcommand{\ccm}{\texttt{cuchebmatrix}}
\newcommand{\ccp}{\texttt{cuchebpoly}}
\newcommand{\ccl}{\texttt{cucheblanczos}}

\def\CC{{C\nolinebreak[4]\hspace{-.05em}\raisebox{.4ex}{\tiny\bf ++}}}

\let\oldding\ding% Make a copy of \ding called \oldding
\renewcommand{\ding}[1]{\ifnum#1=73 $\dag$\else\oldding{#1}\fi}

\begin{document}

\begin{frontmatter}

\title{Cucheb: a GPU implementation of the filtered Lanczos
procedure\tnoteref{mytitlenote}}
\tnotetext[mytitlenote]{This work was partially supported by the European
Research Council under the European Union's Seventh Framework Programme
(FP7/2007–2013)/ERC grant agreement no. 291068. The views expressed in this
article are not those of the ERC or the European Commission, and the European
Union is not liable for any use that may be made of the information contained
here; and by the Scientific Discovery through Advanced Computing (SciDAC)
program funded by U.S. Department of Energy, Office of Science, Advanced
Scientific Computing Research and Basic Energy Sciences DE-SC0008877.}

\author[a]{Jared L. Aurentz}
\author[b]{Vassilis Kalantzis\corref{author}}
\author[b]{Yousef Saad}

\cortext[author] {Corresponding author.\\\textit{E-mail address:}
kalantzi@cs.umn.edu}
\address[a]{Mathematical Institute, University of Oxford, Andrew Wiles Building,
Woodstock Road, OX2 6GG, Oxford, UK.}
\address[b]{Department of Computer Science and Engineering, University of
Minnesota, 4-192 EE/CS Bldg., 200 Union Street S.E., Minneapolis, Minnesota,
55455, USA.}

%%%%%%%%%%%%%%%%%%%%%%%%%%%%%%%%%%%%%%%%%%%%%%%%%%%%%%%%%%%%%%%%%%%%%%%%%%%%%%%%
% Abstract
%%%%%%%%%%%%%%%%%%%%%%%%%%%%%%%%%%%%%%%%%%%%%%%%%%%%%%%%%%%%%%%%%%%%%%%%%%%%%%%%
\begin{abstract}
This paper describes the software package Cucheb, a GPU implementation of the
filtered Lanczos procedure for the solution of large sparse symmetric
eigenvalue problems. The filtered Lanczos procedure uses a carefully chosen
polynomial spectral transformation to accelerate convergence of the Lanczos
method when computing eigenvalues within a desired interval. This method has
proven particularly effective for eigenvalue problems that arise in electronic
structure calculations and density functional theory. We compare our
implementation against an equivalent CPU implementation and show that using the
GPU can reduce the computation time by more than a factor of 10.
\end{abstract}

%%%%%%%%%%%%%%%%%%%%%%%%%%%%%%%%%%%%%%%%%%%%%%%%%%%%%%%%%%%%%%%%%%%%%%%%%%%%%%%%
% Keywords
%%%%%%%%%%%%%%%%%%%%%%%%%%%%%%%%%%%%%%%%%%%%%%%%%%%%%%%%%%%%%%%%%%%%%%%%%%%%%%%%

\end{frontmatter}

%\linenumbers
%% end of elsevier

%%%%%%%%%%%%%%%%%%%%%%%%%%%%%%%%%%%%%%%%%%%%%%%%%%%%%%%%%%%%%%%%%%%%%%%%%%%%%%%%
% AMS Subject Classification
%%%%%%%%%%%%%%%%%%%%%%%%%%%%%%%%%%%%%%%%%%%%%%%%%%%%%%%%%%%%%%%%%%%%%%%%%%%%%%%%
%\begin{AMS}
%65F15, 65F50, 65Y05, 68W10 
%\end{AMS}

%%%%%%%%%%%%%%%%%%%%%%%%%%%%%%%%%%%%%%%%%%%%%%%%%%%%%%%%%%%%%%%%%%%%%%%%%%%%%%%%
% Headings
%%%%%%%%%%%%%%%%%%%%%%%%%%%%%%%%%%%%%%%%%%%%%%%%%%%%%%%%%%%%%%%%%%%%%%%%%%%%%%%%
\pagestyle{myheadings} \thispagestyle{plain}%
\markboth{Jared\ L.\ Aurentz, Vassilis\ Kalantzis, Yousef\ Saad}{Cucheb: GPU
accelerated filtered Lanzcos}

\newpage
\begin{small}
\noindent
{\em Program title:} Cucheb                  \\
{\em Licensing provisions:}    MIT            \\
{\em Programming language:} CUDA C/C++         \\
{\em Nature of problem:} Electronic structure calculations require the
computation of all eigenvalue-eigenvector pairs of a symmetric matrix that lie
inside a user-defined real interval.\\
{\em Solution method:} To compute all the eigenvalues within a given interval a
polynomial spectral transformation is constructed that maps the desired
eigenvalues of the original matrix to the exterior of the spectrum of the
transformed matrix. The Lanczos method is then used to compute the desired
eigenvectors of the transformed matrix, which are then used to recover the
desired eigenvalues of the original matrix. The bulk of the operations are
executed in parallel using a graphics processing unit (GPU).\\
{\em Runtime:} Variable, depending on the number of eigenvalues sought and
the size and sparsity of the matrix.
\end{small}

\section{Introduction}

This paper describes the software package \emph{Cucheb}, a GPU implementation
of the filtered Lanczos procedure \cite{FanSaa2012}.  The filtered Lanczos
procedure (FLP) uses carefully chosen polynomial spectral transformations to
accelerate the computation of all the eigenvalues and corresponding
eigenvectors of a real symmetric matrix $A$ inside a given interval. The chosen
polynomial maps the eigenvalues of interest to the extreme part of the spectrum
of the transformed matrix. The Lanczos method \cite{Lan1950} is then applied to
the transformed matrix which typically converges quickly to the invariant
subspace corresponding to the extreme part of the spectrum. This technique has
been particularly effective for large sparse eigenvalue problems arising in
electronic structure calculations \cite{SchCheSaa2012, Zho2010, ZhoSaa2007,
ZhoSaaTiaChe2006, SaaStaCheWuOgu1996}. 

In the density functional theory framework (DFT) the solution of the
all-electron Schr\"odinger equation is replaced by a one-electron Schr\"odinger
equation with an effective potential which leads to a nonlinear eigenvalue
problem known as the Kohn-Sham equation \cite{PhysRev.136.B864,
PhysRev.140.A1133}:
\begin{equation}\label{eq:kseq}
\left[ 
  - \frac{\nabla^2}{2} + V_{ion}(r) + V_H(\rho(r), r) + V_{XC}(\rho(r), r)  \right]
     \Psi_i(r) = E_i \Psi_i(r),
\end{equation}
where $\Psi_i(r)$ is a wave function and $E_i$ is a Kohn-Sham eigenvalue. The
ionic potential $V_{ion}$ reflects contributions from the core and depends on
the position $r$ only. Both the Hartree and the exchange-correlation
potentials depend on the charge density: 
\begin{equation}\label{eq:rho}
\rho(r) = 2 \sum_{i=1}^{n_{occ}} |\Psi_i(r)|^2,
\end{equation} 
where $n_{occ}$ is the number of occupied states (for most systems of interest
this is half the number of valence electrons). Since the total potential
$V_{total} = V_{ion}+V_H + V_{XC}$ depends on $\rho(r)$ which itself depends on
eigenfunctions of the Hamiltonian, Equation \nref{eq:kseq} can be viewed as a
nonlinear eigenvalue problem or a \emph{nonlinear eigenvector problem}. The
Hartree potential $V_H$ is obtained from $\rho$ by solving the Poisson equation
$ \nabla^2 V_H(r) = -4 \pi \rho(r)$ with appropriate boundary conditions. The
exchange-correlation term $V_{XC}$ is the key to the DFT approach and it
captures the effects of reducing the problem from many particles to a
one-electron problem, i.e., from replacing wavefunctions with many coordinates
into ones that depend solely on space location $r$. 

Self-consistent iterations for solving the Kohn-Sham equation start with an
initial guess of the charge density $\rho(r)$, from which a guess for
$V_{total}$ is computed. Then \nref{eq:kseq} is solved for $\Psi_i(r)$'s and a
new $\rho(r)$ is obtained from \nref{eq:rho} and the potentials are updated.
Then \nref{eq:kseq} is solved again for a new $\rho$ obtained from the new
$\Psi_i(r)$'s, and the process is repeated until the total potential has
converged. 

A typical electronic structure calculation with many atoms requires the
calculation of a large number of eigenvalues, specifically the $n_{occ}$
leftmost ones. In addition, calculations based on time-dependent density
functional theory~\cite{TDDFT-Book,russ-tdlda}, require a substantial number of
unoccupied states, states beyond the Fermi level, in addition to the occupied
ones. Thus, it is not uncommon to see eigenvalue problems in the size of
millions where tens of thousands of eigenvalues may be needed.

Efficient numerical methods that can be easily parallelized in current
high-performance computing environments are therefore essential in electronic
structure calculations. The high computational power offered by GPUs has
increased their presence in the numerical linear algebra community and they are
gradually becoming an important tool of scientific codes for solving
large-scale, computationally intensive eigenvalue problems. While GPUs are
mostly known for their high speedups relative to CPU-bound
operations\footnote{See also the MAGMA project at
http://icl.cs.utk.edu/magma/index.html}, sparse eigenvalue computations can
also benefit from hybrid CPU-GPU architectures. Although published literature
and scientific codes for the solution of sparse eigenvalue problems on a GPU
have not been as common as those that exist for multi-CPU environments, recent
studies conducted independently by some of the authors of this paper
demonstrated that the combination of polynomial filtering eigenvalue solvers
with GPUs can be beneficial \cite{Kal2015, Aur2014}. 

The goal of this paper is twofold. First we describe our open source software
package Cucheb\footnote{\url{https://github.com/jaurentz/cucheb}} that uses the
filtered Lanczos procedure to accelerate large sparse eigenvalue computations
using Nvidia brand GPUs. Then we demonstrate the effectiveness of using GPUs to
accelerate the filtered Lanczos procedure by solving a set of eigenvalue
problems originating from electronic stucture calculations with Cucheb and
comparing it with a similar CPU implementation.

The paper is organized as follows. Section~\ref{sec:flp} introduces the concept of
polynomial filtering for symmetric eigenvalue problems and provides the basic
formulation of the filters used. Section~\ref{sec:implementation} discusses the
proposed GPU implementation of the filtered Lanczos procedure.
Section~\ref{sec:numex} presents computational results with the proposed GPU
implementations. Finally, concluding remarks are presented in
Section~\ref{sec:conclusions}.

\section{The filtered Lanczos procedure} \label{sec:flp}

The Lanczos algorithm and its variants \cite{Lan1950, WuSim2000, BagCalRei2003,
CalReiSor1994, CulDon1974, Pai1971, Sor1996} are well-established methods for
computing a subset of the spectrum of a real symmetric matrix. These methods
are especially adept at approximating eigenvalues lying at the extreme part of
the spectrum \cite{BeaEmbSor2005, Leh1995, Saa1980, Sim1984a}. When the desired
eigenvalues are well inside the spectral interval these techniques can become
ineffective and lead to large computational and memory costs. Traditionally,
this is overcome by mapping interior eigenvalues to the exterior part using a
shift-and-invert spectral transformation (see for example \cite{Saa2011} or
\cite{Wat2007}). While shift-and-invert techniques typically work very well,
they require solving linear systems involving large sparse matrices which can
be difficult or even infeasible for certain classes of matrices.

The filtered Lanczos procedure (FLP) offers an appealing alternative for such
cases. In this approach interior eigenvalues are mapped to the exterior of the
spectrum using a polynomial spectral transformation.  Just as with
shift-and-invert, the Lanczos method is then applied to the transformed matrix
\cite{FanSaa2012}. The key difference is that polynomial spectral
transformations only require matrix-vector multiplication, a task that is often
easy to parallelize for sparse matrices. For FLP constructing a good polynomial
spectral transformation is the most important prepocessing step. 

\subsection{Polynomial spectral transformations}
Let $A \in \bR{n \times n}$ be symmetric and let 
\begin{equation}
A = V\Lambda V^T 
\end{equation}
be its spectral decomposition, where $V \in \bR{n \times n}$ is an orthogonal
matrix and $\Lambda = \textrm{diag}\left(\lambda_1,\ldots,\lambda_n\right)$ is
real and diagonal. A \emph{spectral transformation} of $A$ is a mapping of the
form 
\begin{equation}
f(A) = V f(\Lambda) V^T,
\end{equation}
where $f(\Lambda) = \textrm{diag}\left(f(\lambda_1),\ldots,f(\lambda_n)\right)$
and $f$ is any (real or complex) function $f$ defined on the spectrum of $A$.
Standard examples in eigenvalue computations include the shift-and-invert
transformation $f(z) = (z-\rho)^{-1}$ and $f(z) = z^k$ for subspace iteration. 

A \emph{polynomial spectral transformation} or \emph{filter polynomial} is any
spectral transformation that is also a polynomial. For the filtered Lanczos
procedure a well constructed filter polynomial means rapid convergence and a
good filter polynomial $p$ should satisfy the following requirements: a) the
desired eigenvalues of $A$ are the largest in magnitude eigenvalues of $p(A)$,
b) the construction of $p$ requires minimal knowledge of the spectrum of $A$,
and c) multiplying a vector by $p(A)$ is relatively inexpensive and easy to
parallelize.

Our implementation of the FLP constructs polynomials that satisfy the above
requirements using techniques from digital filter design. The basic idea is to
construct a polynomial filter by approximating an ``ideal'' filter which
maps the desired eigenvalues of $A$ to eigenvalues of largest magnitude in
$p(A)$.

\subsection{Constructing polynomial transformations}
Throughout this section it is assumed that the spectrum of $A$ is contained
entirely in the interval $[-1,1]$. In practice, this assumption poses no
restrictions since the eigenvalues of $A$ located inside the interval
$[\lambda_{\min},\lambda_{\max}]$, where $\lambda_{\min},\ \lambda_{\max}$
denote the algebraically smallest and largest eigenvalues of $A$ respectively,
can be mapped to the interval $[-1,1]$ by the following linear transformation:
\begin{equation}
A := (A - cI)/e,\ \ c =
\dfrac{\lambda_{\min}+\lambda_{\max}}{2},\ \
e=\dfrac{\lambda_{\max}-\lambda_{\min}}{2}. 
\label{transf}
\end{equation} 
Since $\lambda_{\min}$ and $\lambda_{\max}$ are exterior eigenvalues of $A$,
one can obtain very good estimates by performing a few Lanczos steps. We will
see in Section~\ref{sec:numex} that computing such estimates constitutes only a
modest fraction of the total compute time. 

Given a subinterval $\ab \subset [-1,1]$ we wish to compute all eigenvalues of
$A$ in $\ab$ along with their corresponding eigenvectors. Consider first the
following spectral transformation:
\begin{equation}
\phi(z) = 
\left\lbrace
\begin{array}{ccc}
1, & z \in \ab,\\
0, & \mathrm{otherwise}.
\end{array}
\right.
\end{equation}
The function $\phi$ is just an indicator function, taking the value 1 inside
the interval $\ab$ and zero outside. When acting on $A$, $\phi$ maps the
desired eigenvalues of $A$ to the repeated eigenvalue $1$ for $\phi(A)$ and all
the unwanted eigenvalues to $0$. Moreover, the invariant subspace which
corresponds to eigenvalues of $A$ within the interval $\ab$ is identical to the
invariant subspace of $\phi(A)$ which corresponds to the multiple eigenvalue 1.
Thus, applying Lanczos on $\phi(A)$ computes the same invariant subspace, with
the key difference being that the eigenvalues of interest (mapped to one) are
well-separated from the unwanted ones (mapped to zero), and rapid convergence
can be established. Unfortunately, such a transformation is not practically
significant as there is no cost-effective way to multiply a vector by
$\phi(A)$.

A practical alternative is to replace $\phi$ with a polynomial $p$ such that
$p(z) \approx \phi(z)$ for all $z \in [-1,1]$. Such a $p$ will then map the
desired eigenvalues of $A$ to a neighborhood of $1$ for $p(A)$. Moreover, since
$p$ is a polynomial, applying $p(A)$ to a vector only requires matrix-vector
multiplication with $A$. 

In order to quickly construct a $p$ that is a good approximation to $\phi$ it
is important that we choose a good basis. For functions supported on $[-1,1]$
the obvious choice is Chebyshev polynomials of the first kind. Such
representations have already been used successfully for constructing polynomial
spectral transformations and for approximating matrix-valued functions in
quantum mechanics (see for example \cite{Aur2014, BekKokSaa2008, FanSaa2012,
JayKimSaaChe1999, Saa2006, Zho2010, ZhoSaa2007, SchCheSaa2012,
SilRoeVotKre1996, WeiWelAlvFeh2006, ZhoSaaTiaChe2006}). 

\begin{figure}[!h]
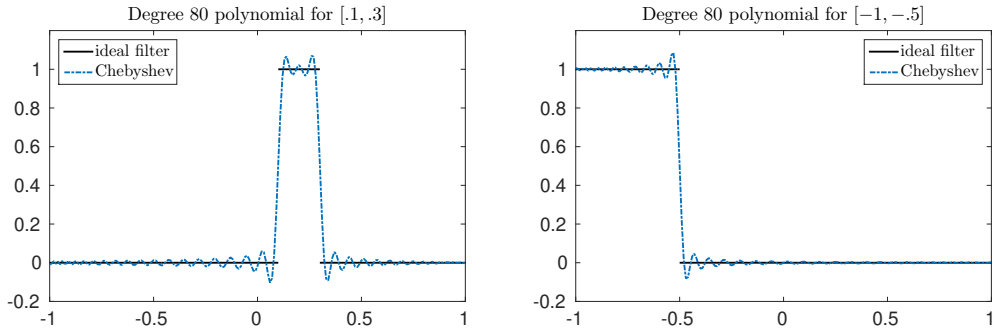
 
\centering
\includegraphics[width=0.45\textwidth]{indicator1a.eps}\hspace{.05\textwidth}
\includegraphics[width=0.45\textwidth]{indicator1b.eps}
\caption{Chebyshev approximation of the ideal filter
$\phi$ using a degree $80$ polynomial. Left: $\ab=[.1,.3]$, right:
$\ab=[-1,-.5]$.}
\label{fig:filters}
\end{figure}

Recall that the Chebyshev polynomials of the first kind obey the following
three-term recurrence
\begin{equation}\label{def:cheb1} 
T_{i+1}(z) = 2zT_i(z) - T_{i-1}(z), \ i \geq 1.
\end{equation} 
starting with $T_0(z) = 1$, $T_1(z) = z$. The Chebyshev polynomials also satisfy
the following orthogonality condition and form a complete orthogonal
set for the Hilbert space $L^2_{\mu}
\left([-1,1]\right)$, $d\mu(z) = \left(1-z^2\right)^{-1/2}dz$:
\begin{equation}
\int_{-1}^1 \frac{T_i(z) T_j(z)}{\sqrt{1-z^2}} dz= 
\left\lbrace
\begin{array}{ccc}
\pi, & i = j = 0,\\
\frac{\pi}{2}, & i = j > 0,\\
0, & \mathrm{otherwise}.
\end{array}
\right.
\end{equation}
Since $\phi \in L^2_{\mu}\left([-1,1]\right)$ it possesses a convergent
Chebyshev series
\begin{equation}
\phi(z) = \sum_{i=0}^{\infty} b_i T_i(z),
\end{equation}
where the $\{b_i\}_{i=0}^{\infty}$ are defined as follows:
\begin{equation}
b_i 
= \frac{2 - \delta_{i0}}{\pi}\int_{-1}^1 \frac{\phi(z) T_i(z)}{\sqrt{1-z^2}} dz,
\end{equation}
where $\delta_{ij}$ represents the Dirac delta symbol. For a given $\alpha$
and $\beta$ the $\{b_i\}$ are known analytically (see for example
\cite{Jac1930}),
\begin{equation}
b_i = 
\left\lbrace
\begin{array}{ccc}
\left(\arccos(\alpha) - \arccos(\beta)\right)/\pi, & i = 0,\\
2\left(\sin\left(i\arccos(\alpha)\right) -
\sin\left(i\arccos(\beta)\right)\right)/i\pi, & i > 0.
\end{array}
\right.
\end{equation}
An obvious choice for constructing $p$ is to fix a degree $m$ and truncate the
Chebyshev series of $\phi$,
\begin{equation}
p_m(z) = \sum_{i=0}^{m} b_i T_i(z).
\end{equation}
Due to the discontinuities of $\phi$, $p_m$ does not converge to $\phi$
uniformly as $m \rightarrow \infty$. The lack of uniform convergence is not an
issue as long as the filter polynomial separates the wanted and unwanted
eigenvalues. Figure~\ref{fig:filters} illustrates two polynomial spectral
transformations constructed by approximating $\phi$ on two different intervals.
Even with the rapid oscillations near the ends of the subinterval, these
polynomials are still good candidates for separating the spectrum.

Figure \ref{fig:filters} shows approximations of the ideal filter $\phi$ for
two different subintervals of $[-1,1]$, using a fixed degree $m=80$. In the
left subfigure the interval of interest is located around the middle of the
spectrum $\ab = [.1,.3]$, while in the right subfigure the interval of interest
is located at the left extreme part $\ab = [-1,-.5]$. Note that the
oscillations near the discontinuities do not prevent the polynomials from
separating the spectrum.

Since $A$ is sparse, multiplying $p(A)$ by a vector can be done efficiently in
parallel using a vectorized version of Clenshaw's algorithm \cite{Cle1955} when
$p$ is represented in a Chebyshev basis. Moreover Clenshaw's algorithm can be run
entirely in real arithmetic whenever the Chebyshev coefficients of $p$ are
real.

\subsection{Filtered Lanczos as an algorithm}

Assuming we've constructed a polynomial filter $p$, we can approximate
eigenvalues of $A$ by first approximating eigenvalues and eigenvectors of
$p(A)$ using a simple version of the Lanczos method \cite{Lan1950}. Many of the
matrices arising in practical applications possess repeated eigenvalues,
requiring the use of block Lanczos algorithm \cite{CulDon1974}, so we describe
the block version of FLP as it contains the standard algorithm as a special
case.

Given a block size $r$ and a matrix $Q \in \bR{n \times r}$ with orthonormal
columns, the filtered Lanzos procedure iteratively constructs an orthonormal
basis for the Krylov subspace generated by $p(A)$ and $Q$: 
\begin{equation}\label{eq:krylov}
\mathcal{K}_k(p(A),Q) = \mathrm{span}\{Q,p(A)Q,\ldots,p(A)^{k-1}Q\}.
\end{equation}

Let us denote by $Q_k \in \bR{n \times rk}$ the matrix whose columns are
generated by $k-1$ steps of the block Lanczos algorithm. Then, for each integer
$k$ we have $Q_k^TQ_k = I$ and
$\mathrm{range}(Q_k)=\mathrm{span}(\mathcal{K}_k(p(A),Q))$. Since $p(A)$ is
symmetric the columns of $Q_k$ can be generated using short recurrences. This
implies that there exists symmetric $\{D_i\}_{i=1}^{k}$ and upper-triangular
$\{S_i\}_{i=1}^{k}$, $D_i, S_i \in \bR{r \times r}$, such that
\begin{equation}\label{eq:lanczos}
p(A)Q_k = Q_{k+1} \tilde{T}_{k},
\end{equation}
where
\begin{equation}
\tilde{T}_{k} = 
\begin{bmatrix}
 T_k \\
 S_k E_k^T    \\
\end{bmatrix},\ \ 
 T_{k} = 
 \begin{bmatrix}
 D_1 & S_1^T & \phantom{\ddots} &  & \\
 S_1 & D_2 & S_2^T & & \phantom{\ddots} \\
 & S_2 & D_3 & \ddots & \\
 &  & \ddots & \ddots & S_{k-1}^T \\
 \phantom{\ddots} &  &  & S_{k-1} & D_k
 \end{bmatrix},
\end{equation}
and $E_k \in \mathbb{R}^{kr\times r}$ denotes the last $r$ columns of the 
identity matrix of size $kr \times kr$. Left multiplying \eqref{eq:lanczos} 
by $Q_k^T$ gives the Rayleigh-Ritz projection
\begin{equation}
Q^T_{k}p(A)Q_k = T_{k}.
\end{equation}
The matrix $T_k$ is symmetric and banded, with a semi-bandwidth of size $r$.
The eigenvalues of $T_k$ are the Ritz values of $p(A)$ associated with the
subspace spanned by the columns of $Q_k$ and for sufficiently large $k$ the
dominant eigenvalues of $p(A)$ will be well approximated by these Ritz values.
Of course we aren't actually interested in the eigenvalues of $p(A)$ but those
of $A$. We can recover these eigenvalues by using the fact that $p(A)$ has the
same eigenvectors as $A$. Assuming that an eigenvector $v$ of $p(A)$ has been
computed accurately we can recover the corresponding eigenvalue $\lambda$ of
$A$ from the Rayleigh quotient of $v$:
\begin{equation}\label{eq:rayleigh}
\lambda = \frac{v^T A v}{v^T v}.
\end{equation}

In practice we will often have only a good approximation $\hat{v}$ of $v$. The
approximate eigenvector $\hat{v}$ will be a Ritz vector of $p(A)$ associated
with $Q_k$. To compute these Ritz vectors we first compute an
eigendecomposition of $T_k$. Since $T_k$ is real and symmetric there exists an
orthogonal matrix $W_k \in \bR{rk \times rk}$ and a diagonal matrix $\Lambda_k
\in \bR{rk \times rk}$ such that
\begin{equation}\label{eq:proj}
T_k W_k = W_k \Lambda_k.
\end{equation}
Combining \eqref{eq:lanczos} and \eqref{eq:proj}, the Ritz vectors of $p(A)$
are formed as $\hat{V}_k = Q_k W_k$.

\section{Cucheb: a GPU implementation of the filtered Lanczos procedure}
\label{sec:implementation}

A key advantage of the filtered Lanczos procedure is that it requires only
matrix-vector multiplication, an operation that uses relatively low memory and
that is typically easy to parallelize compared to solving large linear systems.
FLP and related methods have already been successfully implemented on
multi-core CPUs and distributed memory machines \cite{SchCheSaa2012}. 

\subsection{The GPU architecture}

A graphical processing unit (GPU) is a single instruction multiple data (SIMD)
scalable model which consists of multi-threaded streaming Multiprocessors
(SMs), each one equipped with multiple scalar processor cores (SPs), with each
SP performing the same instruction on its local portion of data. While they
were initially developed for the purposes of graphics processing, GPUs were
adapted in recent years for general purpose computing. The development of the
Compute Unified Device Architecture (CUDA) \cite{NvidiaCuda} parallel
programming model by Nvidia, an extension of the C language, provides an easy
way for computational scientists to take advantage of the GPU's raw power. 

Although the CUDA programming language allows low level access to Nvidia GPUs,
the Cucheb library accesses the GPU through high level routines included as 
part of the Nvidia CUDA Toolkit. The main advantage of this is that one only
has to update to the latest version of the Nvidia's toolkit in order to make
use of the latest GPU technology.

\subsection{Implementation details of the Cucheb software package}

In this section we discuss the details of our GPU implementation of FLP. Our
implementation will consist of a high-level, open source \CC\ library called
\emph{Cucheb} \cite{AurKal2015} which depends only on the Nvidia CUDA Toolkit
\cite{NicBucGarSka2008, NvidiaCuda} and standard \CC\ libraries, allowing for
easy interface with Nvidia brand GPUs. At the user level, the Cucheb software
library consists of three basic data structures:
\begin{itemize}
\item \ccm
\item \ccl
\item \ccp
\end{itemize}
The remainder of this section is devoted to describing the role of each of
these data structures.

\subsubsection{Sparse matrices and the \ccm\ object}

The first data structure, called \ccm, is a container for storing and
manipulating sparse matrices. This data structure consists of two sets 
of pointers, one for data stored in CPU memory and one for data stored 
in GPU memory. Such a duality of data is often necessary for GPU 
computations if one wishes to avoid costly memory transfers between the 
CPU and GPU. To initialize a \ccm\ object one simply passes the path to 
a symmetric matrix stored in the matrix market file format 
\cite{BoiPozRemBarDon1997}. The following segment of Cucheb code illustrates 
how to initialize a \ccm\ object using the matrix \texttt{H2O} downloaded 
from the University of Florida sparse matrix collection
\cite{DavHu2011}:
\begin{verbatim}

#include "cucheb.h"

int main(){

  // declare cuchebmatrix variable
  cuchebmatrix ccm;

  // create string with matrix market file name
  string mtxfile("H2O.mtx");

  // initialize ccm using matrix market file
  cuchebmatrix_init(&mtxfile, &ccm);

  .
  .
  .

}

\end{verbatim}
The function \verb|cuchebmatrix_init| opens the data file, checks that the
matrix is real and symmetric, allocates the required memory on the CPU and GPU,
reads the data into CPU memory, converts it to an appropriate format for the
GPU and finally copies the data into GPU memory. By appropriate format we mean
that the matrix is stored on the GPU in compressed sparse row (CSR) format with
no attempt to exploit the symmetry of the matrix. CSR is used as it is one of
the most generic storage scheme for performing sparse matrix-vector
multiplications using the GPU. (See \cite{RegGil2012,Bel2009} and references
therein for a discussion on the performance of sparse matrix-vector
multiplications in the CSR and other formats.) Once a \ccm\ object has been
created, sparse matrix-vector multiplications can then be performed on the GPU
using the Nvidia CUSPARSE library \cite{NvidiaCusparse}.

\subsubsection{Lanczos and the \ccl\ object}

The second data structure, called \ccl, is a container for storing and
manipulating the vectors and matrices associated with the Lanczos process. As
with the \ccm\ objects, a \ccl\ object possesses pointers to both CPU and GPU
memory. While there is a function for initializing a \ccl\ object, the average
user should never do this explicitly. Instead they should call a higher level
routine like\\ \verb|cuchebmatrix_lanczos| which takes as an argument an
uninitialized \ccl\ object. Such a routine will then calculate an appropriate
number of Lanczos vectors based on the input matrix and initialize the \ccl\
object accordingly. 

Once a \ccm\ object and corresponding \ccl\ object have been initialized, one of
the core Lanczos algorithms can be called to iteratively construct the Lanczos
vectors. Whether iterating with $A$ or $p(A)$, the core Lanczos routines in Cucheb
are essentially the same. The algorithm starts by constructing an orthonormal
set of starting vectors (matrix $Q$ in \eqref{eq:krylov}). Once the vectors are
initialized the algorithm expands the Krylov subspace, peridiocally checking
for convergence. To check convergence the projected problem \eqref{eq:proj} is
copied to the CPU, the Ritz values are computed and the residuals are checked.
If the algorithm has not converged the Krylov subspace is expanded further and
the projected problem is solved again. For stability reasons Cucheb uses full
reorthogonalization to expand the Krylov subspace, making the algorithm more
akin to the Arnoldi method \cite{Arn1951}. Due to the full reorthogonalization,
the projected matrix $T_k$ from \eqref{eq:proj} will not be symmetric exactly
but it will be symmetric to machine precision, which justifies the use of an
efficient symmetric eigensolver (see for example \cite{Par1980}). The cost of
solving the projected problem is negligible compared to expanding the Krylov
subspace, so we can afford to check convergence often. All the operations
required for reorthogonalization are performed on the GPU using the Nvidia
CUBLAS library \cite{NvidiaCublas}. Solving the eigenvalue problem for $T_k$
is done on the CPU using a special purpose built banded symmetric eigensolver 
included in the Cucheb library.

It is possible to use selective reorthogonalization \cite{ParSco1979,
Sim1984a, Sim1984b} or implicit restarts \cite{WuSim2000, BagCalRei2003},
though we don't make use of these techniques in our code. In
Section~\ref{sec:numex} we will see that the dominant cost in the algorithm is
the matrix-vector multiplication with $p(A)$, so reducing the number of
products with $p(A)$ is the easiest way to shorten the computation time.
Techniques like implicit restarting can often increase the number of iterations
if the size of the maximum allowed Krylov subspace is too small, meaning we
would have to perform more matrix-vector multiplications. Our experience
suggests that the best option is to construct a good filter polynomial and then
compute increasingly larger Krylov subspaces until the convergence criterion is
met.

All the Lanczos routines in Cucheb are designed to compute all the eigenvalues in
a user prescribed interval $\ab$. When checking for convergence the Ritz values
and vectors are sorted according to their proximity to $\ab$ and the method is
considered to be converged when all the Ritz values in $\ab$ as well as a few of
the nearest Ritz values outside the interval have sufficiently small residuals.
If the iterations were done using $A$ then the computation is complete and the
information is copied back to the CPU. If the iterations were done with $p(A)$
the Rayleigh quotients are first computed on the GPU and then the information is
copied back to the CPU.

To use Lanczos with $A$ to compute all the eigenvalues in $\ab$ a user is
required to input five variables:
\begin{enumerate}
\item a lower bound on the desired spectrum ($\alpha$)
\item an upper bound on the desired spectrum ($\beta$)
\item a block size
\item an initialized \ccm\ object
\item an uninitialized \ccl\ object
\end{enumerate} 
The following segment of Cucheb code illustrates how to do this using the
function \verb|cuchebmatrix_lanczos| for the interval $\ab = [.5,.6]$, a block
size of $3$ and an already initialized \ccm\ object:
\begin{verbatim}

#include "cucheb.h"

int main(){

  // initialize cuchebmatrix object
  cuchebmatrix ccm;
  string mtxfile("H2O.mtx");
  cuchebmatrix_init(&mtxfile, &ccm);

  // declare cucheblanczos variable
  cucheblanczos ccl;

  // compute eigenvalues in [.5,.6] using block Lanczos 
  cuchebmatrix_lanczos(.5, .6, 3, &ccm, &ccl);

  .
  .
  .

}

\end{verbatim}
This function call will first approximate the upper and lower bounds on the
spectrum of the \ccm\ object. It then uses these bounds to make sure that the
interval $\ab$ is valid. If it is, it will adaptively build up the Krylov
subspace as described above, periodically checking for convergence. For large
matrices or subintervals well inside the spectrum, standard Lanczos may fail to
converge all together. A better choice is to call the routine
\verb|cuchebmatrix_filteredlanczos| which automatically constructs a filter
polynomial and then uses FLP to compute all the eigenvalues in $\ab$.

\subsubsection{Filter polynomials and the \ccp\ object}

To use FLP one needs a way to store and manipulate filter polynomials stored in
a Chebyshev basis. In Cucheb this is done with the \ccp\ object. The \ccp\ object
contains pointers to CPU and GPU memory which can be used to construct and
store filter polynomials. For the filter polynomials from
Section~\ref{sec:flp} one only needs to store the degree, the Chebyshev
coefficients and upper and lower bounds for the spectrum of $A$.

As with \ccl\ objects, a user typically will not need to initialize a \ccp\
object themselves as it will be handled automatically by a higher level
routine. In \verb|cuchebmatrix_filteredlanczos| for example, not only is the
\ccp\ object for the filter polynomial initialized but also the degree at which
the Chebyshev approximation should be truncated is computed. This is done using
a simple formula based on heuristics and verified by experiment. Assuming the
spectrum of $A$ is in $[-1,1]$, a ``good'' degree $m$ for $\ab \subset [-1,1]$
is computed using the following formula:
\begin{equation}\label{eq:optimaldegree}
m = \min\{ m>0: ||p_m-\phi || < \epsilon ||\phi|| \},
\end{equation} 
where $||f||$ is the weighted Chebyshev $2$-norm. The tolerance $\epsilon$ is
a parameter and is chosen experimentally, with the goal of maximizing the
separation power of the filter while keeping the polynomial degree and
consequently the computation time low. 

\begin{figure}[!h]
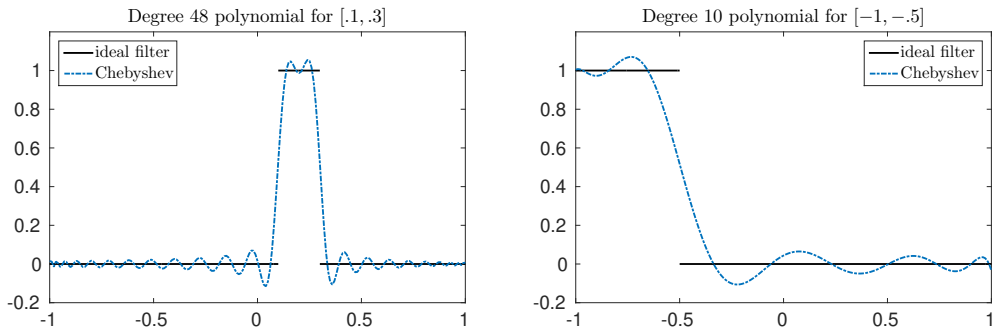
 
\centering
\includegraphics[width=0.45\textwidth]{indicator2a.eps}\hspace{.05\textwidth}
\includegraphics[width=0.45\textwidth]{indicator2b.eps}
\caption{Chebyshev and approximation of the ideal filter
$\phi$. Left: $\ab=[.1,.3]$ with an optimal degree of $48$, right:
$\ab=[-1,-.5]$ with an optimal degree of $10$.}
\label{fig:optimalfilters}
\end{figure}

Figure~\ref{fig:optimalfilters} uses the same ideal filters from
Figure~\ref{fig:filters} but this time computes the filter degree based on
\eqref{eq:optimaldegree}. In the left subfigure the interval of interest is
located around the middle of the spectrum $\ab = [.1,.3]$ and the distance
between $\alpha$ and $\beta$ is relatively small, giving a
filter degree of $48$. In the right subfigure the interval of interest is
located at the left extreme part of the spectrum $\ab = [-1,-.5]$ and the
distance $\alpha$ and $\beta$ is relatively large, giving a
filter degree of $10$. Although these filters seem like worse approximations
than those in Figure~\ref{fig:filters}, the lower degrees lead to much shorter
computation times.

The following segment of Cucheb code illustrates how to use the function\\
\verb|cuchebmatrix_filteredlanczos| to compute all the eigenvalues in the
interval $\ab = [.5,.6]$ of an already initialized \ccm\ object using FLP with
a block size of $3$:
\begin{verbatim}

#include "cucheb.h"

int main(){

  // initialize cuchebmatrix object
  cuchebmatrix ccm;
  string mtxfile("H2O.mtx");
  cuchebmatrix_init(&mtxfile, &ccm);

  // declare cucheblanczos variable
  cucheblanczos ccl;

  // compute eigenvalues in [.5,.6] using filtered Lanczos 
  cuchebmatrix_filteredlanczos(.5, .6, 3, &ccm, &ccl);

  .
  .
  .

}

\end{verbatim}

\section{Experiments} \label{sec:numex}

\begin{figure} 
  \centering
  \includegraphics[width=0.47\textwidth]{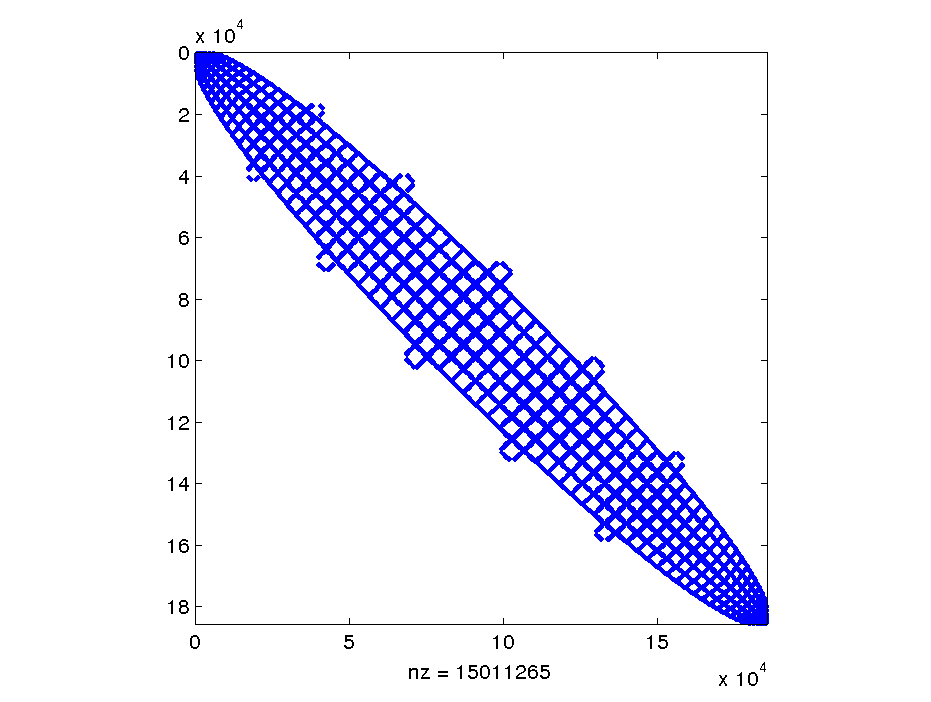}
  \includegraphics[width=0.47\textwidth]{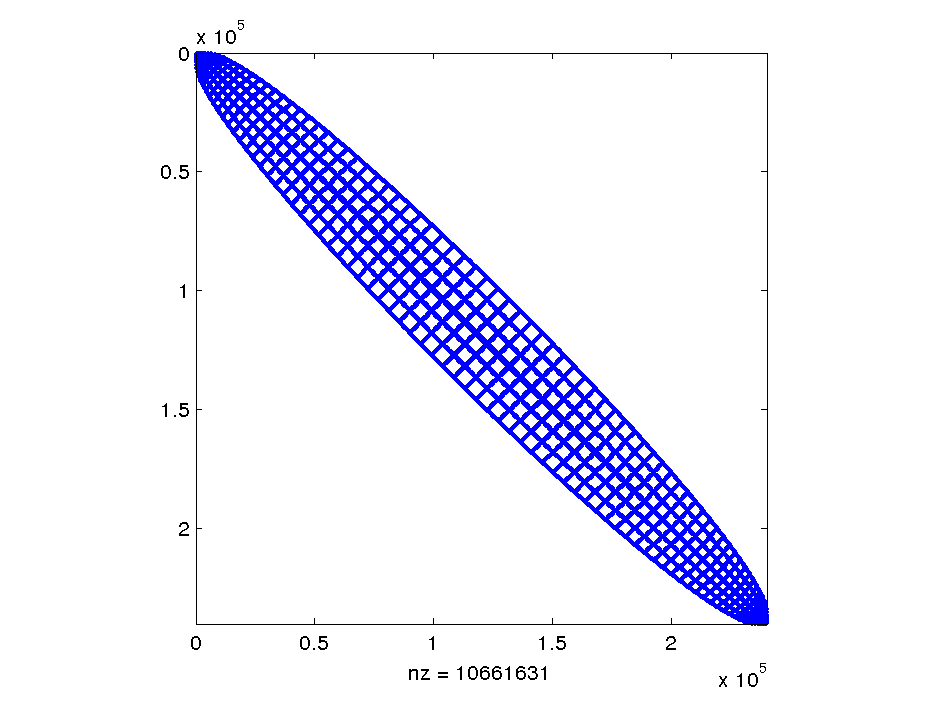}
  \caption{Sparsity pattern of the PARSEC matrices. Left: Si41Ge41H72. Right:
Si87H76.}
  \label{parsec_spy}
\end{figure}

In this section we illustrate the performance of our GPU implementation of the
filtered Lanczos procedure. Our test matrices (Hamiltonians) originate from
electronic structure calculations. In this setting, one is typically interested
in computing a few eigenvalues around the Fermi level of each Hamiltonian. The
Hamiltonians were generated using the PARSEC package \cite{Kroetal2006} and
can be also found in the University of Florida sparse matrix collection
\cite{DavHu2011}.\footnote{https://www.cise.ufl.edu/research/sparse/matrices/}
These Hamiltonians are real, symmetric, and have clustered, as well as multiple,
eigenvalues. Table \ref{tab:matrices} lists the size $n$, the total number of
non-zero entries $nnz$, as well as the endpoints of the spectrum of each
matrix, i.e., the interval defined by the algebraically smallest/largest
eigenvalues. The average number of nonzero entries per row for each Hamiltonian
is quite large, a consequence of the high-order discretization and the addition
of a (dense) `non-local' term. Figure \ref{parsec_spy} plots the sparsity
pattern of matrices Si41Ge41H72 (left) and Si87H76 (right).

All GPU experiments in this section were implemented using the Cucheb library
and performed on the same machine which has an Intel Xeon E5-2680 v3 2.50GHz
processor with 128GB of CPU RAM and two Nvidia K40 GPUs each with 12GB of GPU
RAM and 2880 compute cores. We make no attempt to access mutliple GPUs and all
the experiments were performed using a single K40. 

\begin{table}[!ht]
  \centering
  \begin{tabular}{l|c|c|c|c}
\hline
\multirow{2}{*}{Matrix} & \multirow{2}{*}{$n$} & \multirow{2}{*}{$nnz$} & \multirow{2}{*}{$nnz/n$} & \multirow{2}{*}{Spectral interval} \\
 & & & & \\\hline
\hline
\verb|Ge87H76| & $\phantom{0,{}}112,985$ & $\phantom{0}7,892,195$ & $69.9$ & $[-1.21e{+0},\phantom{-{}}3.28e{+1}]$ \\
\verb|Ge99H100| & $\phantom{0,{}}112,985$ & $\phantom{0}8,451,395$ & $74.8$ & $[-1.23e{+0},\phantom{-{}}3.27e{+1}]$ \\
\verb|Si41Ge41H72| & $\phantom{0,{}}185,639$ & $15,011,265$ & $80.9$ & $[-1.21e{+0},\phantom{-{}}4.98e{+1}]$ \\
\verb|Si87H76| & $\phantom{0,{}}240,369$ & $10,661,631$ & $44.4$ & $[-1.20e{+0},\phantom{-{}}4.31e{+1}]$ \\
\verb|Ga41As41H72| & $\phantom{0,{}}268,096$ & $18,488,476$ & $69.0$ & $[-1.25e{+0},\phantom{-{}}1.30e{+3}]$ \\\hline
\hline
\end{tabular}
  \vspace{2ex}
  \caption{A list of the PARSEC matrices used to evaluate our GPU
implementation, where $n$ is the dimension of the matrix, $nnz$ is the number
of nonzero entries and $[\lambda_{\min},\lambda_{\max}]$ is the spectral
interval.}
  \label{tab:matrices}
\end{table}

Exploiting eigenvalue solvers that are based on matrix factorizations, e.g.,
shift-and-invert techniques, has been shown to be impractical for matrices of
the PARSEC matrix collection \cite{DDFEAST,KalantzisSaa15:SpecSch-TR}. The
reason is that performing the LU factorization of each Hamiltonian results in a
huge amount of fill-in in the associated triangular factors, requiring an
excessive amount of memory and computations \cite{DDFEAST}. On the other hand,
polynomial filtering accesses the Hamiltonians in their original form and only
requires an efficient matrix-vector multiplication routine. Polynomial
filtering has often been reported to be the most efficient numerical method for
solving eigenvalue problems with the PARSEC matrix collection
\cite{SchCheSaa2012, FanSaa2012, Zho2010, ZhoSaa2007, ZhoSaaTiaChe2006,
SaaStaCheWuOgu1996}. This observation led to the development of FILTLAN, a
C/C++ software package which implements the filtered Lanczos procedure with
partial reorthgonalization \cite{FanSaa2012} for serial architectures. The
Cucheb library featured in this paper, although implemented in CUDA, shares
many similarities with FILTLAN. There are, however, a few notable differences.
Cucheb does not implement partial reorthgonalization as is the case in FILTLAN.
Moreover, Cucheb includes the ability to use block counterparts of the Lanczos
method which can be more efficient in the case of multiple or clustered
eigenvalues. Morover FILTLAN uses a more complicated least-squares filter
polynomial while Cucheb utilizes the fitlers described in section~\ref{sec:flp}.

\subsection{GPU benchmarking}

The results of the GPU experiments are summarized in Table~\ref{tab:parsec}. 
The variable `interval' for each
Hamiltonian was set so that it included roughly the same number of eigenvalues
from the left and right side of the Fermi level, and in total 
`eigs' eigenvalues. For each matrix and interval $\ab$ we repeated the same
experiment five times, each time using a different degree $m$ for the filter
polynomial. The variable `iters' shows the number of FLP iterations, while `MV'
shows the total number of matrix-vector products (MV) with $A$, which is
computed using the formula `MV'$\phantom{0}= rm \times \phantom{0}$`iters'.
Throughout this section, the block size of the FLP will be equal to $r=3$.
Finally, the variables `time' and `residual' show the total compute time and
maximum relative residual of the computed eigenpairs. The first four rows for
each matrix correspond to executions where the degree $m$ was selected a
priori. The fifth row corresponds to an execution where the degree was selected
automatically by our implementation, using the mechanism described in
\eqref{eq:optimaldegree}. As expected, using larger values for $m$ leads to
faster convergence in terms of total iterations, since higher degree filters
are better at separating the wanted and unwanted portions of the spectrum.
Although larger degrees lead to less iterations, the amount of work in each
filtered Lanczos iteration is also increasing proportionally. This might lead to
an increase of the actual computational time, an effect verified for each one
of the matrices in Table~\ref{tab:parsec}. 

\begin{table}[!ht]
  \centering
  \resizebox{\textwidth}{!}{%
  \begin{tabular}{l|c|c|c|c|c|c|c}
\hline
\multirow{2}{*}{Matrix} & \multirow{2}{*}{interval} & \multirow{2}{*}{eigs} & \multirow{2}{*}{$m$} & \multirow{2}{*}{iters} & \multirow{2}{*}{MV} & \multirow{2}{*}{time} & \multirow{2}{*}{residual} \\
 & & & & & & & \\\hline
\hline
 & & & $\phantom{0}50$ & $210$ & $\phantom{0}31,500$ & $\phantom{0}31$ & $1.7e{-14}$ \\
 & & & $100$ & $180$ & $\phantom{0}54,000$ & $\phantom{0}40$ & $4.0e{-13}$ \\
\verb|Ge87H76| & $[-0.645,-0.0053]$ & $212$ & $150$ & $150$ & $\phantom{0}67,500$ & $\phantom{0}44$ & $7.4e{-14}$ \\
 & & & $200$ & $150$ & $\phantom{0}90,000$ & $\phantom{0}56$ & $6.3e{-14}$ \\
 & & & $\phantom{0}49$ & $210$ & $\phantom{0}30,870$ & $\phantom{0}31$ & $9.0e{-14}$ \\\hline
 & & & $\phantom{0}50$ & $210$ & $\phantom{0}31,500$ & $\phantom{0}32$ & $6.2e{-13}$ \\
 & & & $100$ & $180$ & $\phantom{0}54,000$ & $\phantom{0}41$ & $8.6e{-13}$ \\
\verb|Ge99H100| & $[-0.650,-0.0096]$ & $250$ & $150$ & $180$ & $\phantom{0}81,000$ & $\phantom{0}56$ & $5.0e{-13}$ \\
 & & & $200$ & $180$ & $108,000$ & $\phantom{0}70$ & $1.1e{-13}$ \\
 & & & $\phantom{0}49$ & $210$ & $\phantom{0}30,870$ & $\phantom{0}32$ & $3.2e{-13}$ \\\hline
 & & & $\phantom{0}50$ & $210$ & $\phantom{0}31,500$ & $\phantom{0}56$ & $6.4e{-13}$ \\
 & & & $100$ & $180$ & $\phantom{0}54,000$ & $\phantom{0}73$ & $2.0e{-11}$ \\
\verb|Si41Ge41H72| & $[-0.640,-0.0028]$ & $218$ & $150$ & $180$ & $\phantom{0}81,000$ & $\phantom{0}99$ & $5.6e{-14}$ \\
 & & & $200$ & $150$ & $\phantom{0}90,000$ & $104$ & $5.0e{-13}$ \\
 & & & $\phantom{0}61$ & $180$ & $\phantom{0}32,940$ & $\phantom{0}52$ & $8.9e{-13}$ \\\hline
 & & & $\phantom{0}50$ & $150$ & $\phantom{0}22,500$ & $\phantom{0}38$ & $3.5e{-14}$ \\
 & & & $100$ & $\phantom{0}90$ & $\phantom{0}27,000$ & $\phantom{0}35$ & $4.0e{-15}$ \\
\verb|Si87H76| & $[-0.660,-0.3300]$ & $107$ & $150$ & $120$ & $\phantom{0}54,000$ & $\phantom{0}63$ & $9.1e{-15}$ \\
 & & & $200$ & $\phantom{0}90$ & $\phantom{0}54,000$ & $\phantom{0}60$ & $1.3e{-13}$ \\
 & & & $\phantom{0}98$ & $\phantom{0}90$ & $\phantom{0}26,460$ & $\phantom{0}35$ & $1.2e{-14}$ \\\hline
 & & & $200$ & $240$ & $144,000$ & $225$ & $1.5e{-15}$ \\
 & & & $300$ & $180$ & $162,000$ & $236$ & $2.1e{-15}$ \\
\verb|Ga41As41H72| & $[-0.640,0.0000]$ & $201$ & $400$ & $180$ & $216,000$ & $306$ & $2.5e{-15}$ \\
 & & & $500$ & $180$ & $270,000$ & $375$ & $1.0e{-12}$ \\
 & & & $308$ & $180$ & $166,320$ & $242$ & $1.5e{-15}$ \\\hline
\end{tabular}}
  \vspace{2ex}
  \caption{Computing the eigenpairs inside an interval using FLP with various
  filter polynomial degrees. Times listed are in seconds.}
  \label{tab:parsec}
\end{table}

Table~\ref{tab:parsec2} compares the percentage of total computation time
required by the different subprocesses of the FLP method. We denote the
preprocessing time, which consists solely of approximating the upper and lower
bounds of the spectrum for $A$, by `PREPROC'. We also denote the total amount
of time spent in the full reorthogonalization and the total amount of time
spent in performing all MV products of the form $p(A)v$ on the GPU, by `ORTH'
and `MV' respectively. As we can verify, all matrices in this experiment
devoted no more than $12$\% of the total compute time to estimating the spectral
interval (i.e. the eigenvalues $\lambda_{\min}$ and $\lambda_{\max}$). For each
one of the PARSEC test matrices, the dominant cost came from the MV products,
due to their relatively large number of non-zero entries. Note that using a
higher degree $m$ will shift the cost more towards the MV products, since the
Lanczos procedure will typically converge in fewer outer steps and thus the
orthogonalization cost reduces.
\begin{table}[!ht]
  \centering
  \begin{tabular}{l|c|c|c|c|c}
\hline
\multirow{2}{*}{Matrix} & \multirow{2}{*}{$m$} & \multirow{2}{*}{iters} & \multirow{2}{*}{PREPROC} & \multirow{2}{*}{ORTH} & \multirow{2}{*}{MV} \\
 & & & & & \\\hline
\hline
  & $\phantom{0}50$ & $210$ & $\phantom{0}7$\% & $22$\% & $52$\%\\
  & $100$ & $180$ & $\phantom{0}5$\% & $13$\% & $71$\%\\
\verb|Ge87H76| & $150$ & $150$ & $\phantom{0}5$\% & $\phantom{0}9$\% & $80$\%\\
  & $200$ & $150$ & $\phantom{0}4$\% & $\phantom{0}7$\% & $84$\%\\
  & $\phantom{0}49$ & $210$ & $\phantom{0}7$\% & $21$\% & $52$\%\\\hline
  & $\phantom{0}50$ & $210$ & $\phantom{0}7$\% & $21$\% & $53$\%\\
  & $100$ & $180$ & $\phantom{0}5$\% & $13$\% & $71$\%\\
\verb|Ge99H100| & $150$ & $180$ & $\phantom{0}4$\% & $10$\% & $79$\%\\
  & $200$ & $180$ & $\phantom{0}3$\% & $\phantom{0}8$\% & $83$\%\\
  & $\phantom{0}49$ & $210$ & $\phantom{0}7$\% & $21$\% & $53$\%\\\hline
  & $\phantom{0}50$ & $210$ & $10$\% & $19$\% & $55$\%\\
  & $100$ & $180$ & $\phantom{0}8$\% & $12$\% & $72$\%\\
\verb|Si41Ge41H72| & $150$ & $180$ & $\phantom{0}6$\% & $\phantom{0}9$\% & $80$\%\\
  & $200$ & $150$ & $\phantom{0}5$\% & $\phantom{0}6$\% & $84$\%\\
  & $\phantom{0}61$ & $180$ & $11$\% & $17$\% & $61$\%\\\hline
  & $\phantom{0}50$ & $150$ & $11$\% & $22$\% & $54$\%\\
  & $100$ & $\phantom{0}90$ & $12$\% & $12$\% & $70$\%\\
\verb|Si87H76| & $150$ & $120$ & $\phantom{0}7$\% & $10$\% & $78$\%\\
  & $200$ & $\phantom{0}90$ & $\phantom{0}7$\% & $\phantom{0}7$\% & $83$\%\\
  & $\phantom{0}98$ & $\phantom{0}90$ & $12$\% & $13$\% & $70$\%\\\hline
  & $200$ & $240$ & $\phantom{0}4$\% & $\phantom{0}8$\% & $82$\%\\
  & $300$ & $180$ & $\phantom{0}4$\% & $\phantom{0}5$\% & $88$\%\\
\verb|Ga41As41H72| & $400$ & $180$ & $\phantom{0}3$\% & $\phantom{0}4$\% & $91$\%\\
  & $500$ & $180$ & $\phantom{0}2$\% & $\phantom{0}3$\% & $93$\%\\
  & $308$ & $180$ & $\phantom{0}4$\% & $\phantom{0}5$\% & $89$\%\\\hline
\end{tabular}
  \vspace{2ex}
  \caption{Percentage of total compute time required by various components of
  the algorithm. For all these examples the dominant computational cost are the
  matrix-vector multiplications (MV).}
  \label{tab:parsec2}
\end{table}

We would like to note that the Cucheb software package is capable of running
Lanczos without filtering. We originally intended to compare filtered Lanczos
with standard Lanczos on the GPU, however for the problems considered in this
paper the number of Lanczos vectors required for convergence exceeded the
memory of the K40 GPU. This suggests that for these particular problems
filtering is not only beneficial for performance but also necessary if this
particular hardware is used.

\subsection{CPU-GPU comparison}
\begin{figure}[!ht] 
  \centering
  \includegraphics[width=0.9\textwidth]{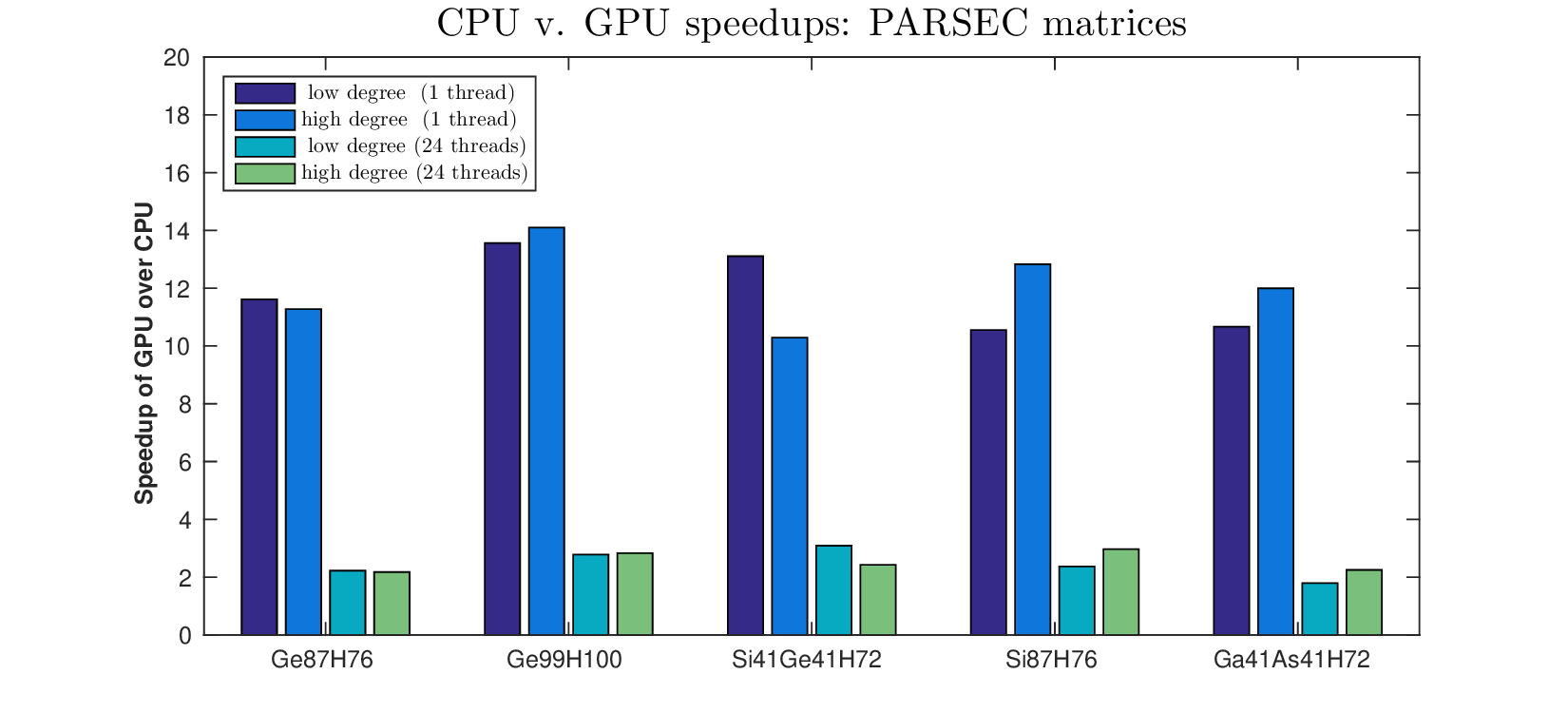}
  \caption{Speedup of the GPU FLP implementation over the CPU (FILTLAN) 
  for the PARSEC test matrices.} 
  \label{fig:gpu_speedup1}
\end{figure}

Figure \ref{fig:gpu_speedup1} shows the speedup of the GPU FLP implementation
over the CPU-based counterpart. The CPU results were obtained by executing the
FILTLAN software package on the {\tt Mesabi} linux
cluster at University of Minnesota Supercomputing Institute. {\tt Mesabi}
consists of 741 nodes of various configurations with a total of 17,784
compute cores provided by Intel Xeon E5-2680 v3 processors. Each node
features two sockets, each socket with twelve physical cores, and each core
with a clock speed of 2.50 GHz. Each node is also equipped with 64 GB of RAM
memory. The FILTLAN package has the option to link the Intel Math Kernel
Library (MKL) when it, as well as a compatible Intel compiler are available. 
For these experiments we used the Intel compiler {\tt icpc} version 11.3.2.

We have divided the comparison into four parts: a ``low degree'' situation when
$m=50$ ($m=200$ for \texttt{Ga41As41H72}), and a ``high degree'' situation when
$m=100$ ($m=300$ for \texttt{Ga41As41H72}) and within each of these we also 
executed FILTLAN using both $1$ thread and $24$ threads. The multithreading was 
handled entirely by the MKL. In the single thread case, the GPU
implementation obtains a speedup which ranges between $10$ and $14$. In the 
$24$ thread 
case, which corresponds to one thread per core on this machine, the speedups 
ranged between $2$ and $3$. 

\section{Conclusion} \label{sec:conclusions}
In this work we presented a GPU implementation of the filtered Lanczos procedure
for solving large and sparse eigenvalue problems such as those that arise from
real-space DFT methods in electronic structure
calculations. Our experiments indicate that the use of GPU architectures in the
context of electronic structure calculations can provide a speedup of at least
a factor of $10$ over a single core CPU implementation and at least of factor of
$2$ for a $24$ core implementation.

Possible future research directions include the utilization of more than one
GPU to perform the filtered Lanczos procedure in computing environments with
access to multiple GPUs. Each GPU can then be used to either perform
the sparse matrix-vector products and other operations of the FLP in parallel,
or compute all eigenpairs in a sub-interval of the original interval. In the later
case the
implementation proposed in this paper can be used without any modifications.
Another interesting extension would be to use additional customization and add
support for other sparse matrix formats. A dense matrix
version of the proposed implementation would also be of interest for solving
sequences of eigenvalue problems as in \cite{BerWorDiN2015}.

\bibliography{intro_paper}

\end{document}